\documentclass[reqno,12pt]{amsart}
\usepackage{url}
\usepackage{amssymb,amsthm,amsmath,amsfonts,amstext,mathrsfs}
\usepackage{latexsym}
\usepackage{epsfig}
\usepackage[all]{xy}
\usepackage{fancyhdr}
\usepackage[unicode]{hyperref}

\hypersetup{colorlinks=true,
linkcolor=blue,
citecolor=blue,
urlcolor=blue,
filecolor=blue,
bookmarksnumbered=true,
pdfstartview=FitH,
pdfhighlight=/N
}

\newtheorem{thm}{Theorem}
\newtheorem{prop}{Proposition}
\newtheorem{lem}{Lemma}
\newtheorem{corr}{Corollary}
\newtheorem{prob}{Problem}

	{\medbreak\par\noindent{{\bf Note.}}}%
	{\hfill{$\Box$}
	\par\medbreak}

\newcommand{\p}{\partial}

\newcommand{\m}{\mathbf}
\newcommand{\mf}{\mathfrak}
\def\d{\,{\rm{d}}}

\title[Giedrius Alkauskas]
{Beltrami vector fields with an icosahedral symmetry} 
\author[Giedrius Alkauskas]{Giedrius Alkauskas}
\address{Vilnius University, Department of Mathematics and Informatics, Naugarduko 24, LT-03225 Vilnius, Lithuania}
\email{giedrius.alkauskas@mif.vu.lt}

\begin{document}
\begin{abstract}A vector field is called a \emph{Beltrami vector field}, if $\m{B}\times(\nabla\times\m{B})=0$. In this paper we construct two unique Beltrami vector fields $\mf{I}$ and $\mf{Y}$, such that $\nabla\times\mf{I}=\mf{I}$, $\nabla\times\mf{Y}=\mf{Y}$, and such that both have an orientation-preserving icosahedral symmetry. Both of them have an additional symmetry with respect to a non-trivial automorphism of the number field $\mathbb{Q}(\,\sqrt{5}\,)$. 
\end{abstract}

\pagestyle{fancy}
\fancyhead{}
\fancyhead[LE]{{\sc Beltrami vector fields}}
\fancyhead[RO]{{\sc G. Alkauskas}}
\fancyhead[CE,CO]{\thepage}
\fancyfoot{}

\date{26 December, 2017}
\thanks{The research of the author was supported by the Research Council of Lithuania grant No. MIP-072/2015}

\maketitle
The aim of this paper is to construct a non-zero solution to the linear system of first order PDE's, given by
\begin{eqnarray*}
\def\arraystretch{2.5}
\left\{\begin{array}{l@{\qquad}l}
\displaystyle{u=\frac{\p w}{\p y}-\frac{\p v}{\p z}},\\
\displaystyle{v=\frac{\p u}{\p z}-\frac{\p w}{\p x}},\\
\displaystyle{w=\frac{\p v}{\p x}-\frac{\p u}{\p y}},
\end{array}\right.
\end{eqnarray*}
such that the solution is in a closed-form, is as simple as possible, is defined in the whole $\mathbb{R}^{3}$, and has an icosahedral symmetry. This has a limited physical applicability, but a huge mathematical interest - see Problems \ref{prob-nt}, \ref{prob-rb}, \ref{prob-is}. For few comments concerning solutions in bounded domains with a boundary condition ensuring that a particle in the field remains inside the domain, see Section \ref{overview}. 
\section{The main results}
In the most general case, if $(M,g)$ is a Riemannian $3$-manifold, the \emph{curl operator} on differential $1$-forms is defined as $*\d:\Omega^{1}(M)\mapsto\Omega^{1}(M)$, where $*$ is the Hodge star operator, and $d$ is the exterior derivative \cite{gadea}. The interpretation in terms of vector fields is as follows. Namely, let $\iota_{X}$ be tensor contraction with respect to a vector field $X$, $d$ be an exterior derivative, and $^{\sharp}$ denote the isomorphism from $1$-forms to vector fields derived from $g$ \cite{gadea,hitchin}. Then we define curl of a vector field $X$, denoted by $\nabla\times X$, by $\nabla\times X=(*\d\iota_{X}g)^{\sharp}$ \cite{etnyre,etnyre2}. For the functional-analytic and spectral treatment of curl operator, see \cite{hiptmair}. In this paper we deal with the most basic example $M=\mathbb{R}^{3}$, and the standard Euclidean metric. But see the end of this paper for questions related to the case $M=\m{S}^{3}$, the three-dimensional sphere. For a vector field $a_{x}\frac{\p}{\p x}+a_{y}\frac{\p}{\p y}+a_{z}\frac{\p}{\p z}$ we usually write $(a_{x},a_{y},a_{z})$.\\

 Thus, let, as usual, $\nabla\times\m{B}=\mathrm{curl}\,\m{B}$. A $3$-dimensional vector field $\m{B}$ in $\mathbb{R}^{3}$ is called a \emph{Beltrami field}, or \emph{force-free magnetic field} (FFF) in physics, if $\m{B}\times(\nabla\times\m{B})=0$; sometimes the condition $\nabla\cdot \m{B}=\mathrm{div}\,\m{B}=0$ is also included. Thus, $\nabla\times\m{B}=f\m{B}$, where $f(x,y,z)$ is a scalar function. Note that for $f$ constant (such fields are called \emph{constant Beltrami fields}), $\nabla\cdot\m{B}=0$ is automatic (here and henceforth we deal with smooth fields only). An important special case of this is when $f\equiv 1$, or $\nabla\times\m{B}=\m{B}$. This is a curl operator eigenvalue $1$ case. Sometimes such vector fields are called \emph{Trkalian vector fields} \cite{baldwin,saygili1,saygili2}. The equality $\nabla\times\m{B}=\m{B}$ is exactly the topic of the current paper. For general smooth $3$-manifolds, Beltrami fields are defined via differential Beltrami $1$-forms \cite{dahl,etnyre}. \\
 
For motivation from a point of view of mathematical physics, we note that in the setting of a magnetohydrodynamic description of plasmas \cite{priest}, force-free magnetic fields play a prominent role. Their characteristic property is the collinearity of magnetic field and electric current and therefore the vanishing of the Lorentz-force. Magnetic field $\m{B}$ and current density $\m{j}$ then satisfy (after normalization)
\begin{eqnarray*}
\left\{\begin{array}{c@{\qquad}l}
\m{j}=f\m{B},\\
\mathrm{curl}\,\m{B}=\m{j},\\
\mathrm{div}\,\m{B}=0.
\end{array}\right.
\end{eqnarray*}
This shows that $\m{B}$ is a Beltrami field. We further quote \cite{tahar}: ``The simultaneous appearence of force-free fields in a great diversity of physical domains has increased their importance in the last decades.  They are particularly a subject of an intensive research in solar physics; indeed, on the sun's surface, most  of  the  observed  structures  and  phenomema (flares, mass ejection, coronal heating, prominences) are mainly due to the strong magnetic field $\m{B}$. When the equilibrium holds - or in a quasi-static evolution - the only significant force, say the Lorentz force $\m{j}\times\m{B}$, must vanish. Therefore,  the electric currents $\m{j}$ are parallel to $\m{B}$ and, thanks to Amp\`{e}re's Law, that means that $\m{B}$ is force-free". And quoting \cite{baldwin}: ``Beltrami fields play a prominent role in the theory of exact, closed form solutions to the Euler and Navier-Stokes equations and their relations to the elctromagnetic wave equations. Moreover, Beltrami fields are related to minimum energy plasma fields and have theorefore garnered much attention from the magnetohydrodynamics community". \\ 

After these motivating remarks, now we will pass the the main topic of this paper. Let $\phi=\frac{1+\sqrt{5}}{2}$. Define 
\begin{eqnarray}
\alpha=\begin{pmatrix}
-1 & 0 & 0\\
0 & -1 & 0\\
0 & 0 & 1 
\end{pmatrix},\quad
\beta=\begin{pmatrix}
0 & 0 & 1\\
1 & 0 &0\\
0 & 1 & 0 
\end{pmatrix},\quad
\gamma=\begin{pmatrix}
\frac{1}{2} & -\frac{\phi}{2} & \frac{1}{2\phi}\\
\frac{\phi}{2} & \frac{1}{2\phi} & -\frac{1}{2}\\
\frac{1}{2\phi} & \frac{1}{2} & \frac{\phi}{2} 
\end{pmatrix}.
\label{gen-ico}
\end{eqnarray}
These three matrices generate the icosahedral group $\mathbb{I}\subset SO(3)$ of order $60$, and $\alpha,\beta$ generate the tetrahedral group $\mathbb{T}$ of order $12$. $\{I,\alpha,\beta^{2}\alpha\beta,\beta\alpha\beta^{2}\}$ form a Klein four group $\mathbb{K}$, which we will encounter later; see (\ref{klein}).\\

Polyhedral symmetries are ubiquitous in science, nature, history, and Art - they manifest from octahedral symmetries of carved stone balls from late Neolithic (c. 3000 BC) found in Scotland, to icosahedral symmetry of capsids of adenoviruses \cite{harnad}. \\

 The first main result of this paper can be stated immediately as follows.
\begin{thm}
Let us define the vector field $\mf{V}=\big{(}\mf{V}_{x},\mf{V}_{y},\mf{V}_{z}\big{)}$, $\mf{V}_{y}=\mf{V}_{x}(y,z,x)$, $\mf{V}_{z}=\mf{V}_{x}(z,x,y)$, where\scriptsize
\begin{eqnarray*}
\mf{V}_{x}&=&2x\sin\Big{(}\frac{x}{2}\Big{)}\sin\Big{(}\frac{\phi y}{2}\Big{)}\sin\Big{(}\frac{z}{2\phi}\Big{)}
-2\phi x\sin\Big{(}\frac{x}{2\phi}\Big{)}\sin\Big{(}\frac{y}{2}\Big{)}\sin\Big{(}\frac{\phi z}{2}\Big{)}+2\phi^{-1} x\sin\Big{(}\frac{\phi x}{2}\Big{)}\sin\Big{(}\frac{y}{2\phi}\Big{)}\sin\Big{(}\frac{z}{2}\Big{)}\\
&+&y\sin z+2y\cos\Big{(}\frac{x}{2}\Big{)}\cos\Big{(}\frac{\phi y}{2}\Big{)}\sin\Big{(}\frac{z}{2\phi}\Big{)}
-2y\cos\Big{(}\frac{x}{2\phi}\Big{)}\cos\Big{(}\frac{y}{2}\Big{)}\sin\Big{(}\frac{\phi z}{2}\Big{)}\\
&+&z\sin y-2z\cos\Big{(}\frac{x}{2}\Big{)}\sin\Big{(}\frac{\phi y}{2}\Big{)}\cos\Big{(}\frac{z}{2\phi}\Big{)}+2z\cos\Big{(}\frac{\phi x}{2}\Big{)}\sin\Big{(}\frac{y}{2\phi}\Big{)}\cos\Big{(}\frac{z}{2}\Big{)},
\end{eqnarray*}
\normalsize
and the vector field $\mf{W}=\big{(}\mf{W}_{x},\mf{W}_{y},\mf{W}_{z}\big{)}$, $\mf{W}_{y}=\mf{W}_{x}(y,z,x)$, $\mf{W}_{z}=\mf{W}_{x}(z,x,y)$, where \scriptsize
\begin{eqnarray*}
\mf{W}_{x}&=&x\cos y-x\cos z\\
&-&\sqrt{5}x\cos\Big{(}\frac{x}{2}\Big{)}\cos\Big{(}\frac{\phi y}{2}\Big{)}\cos\Big{(}\frac{z}{2\phi}\Big{)}+\phi x\cos\Big{(}\frac{x}{2\phi}\Big{)}\cos\Big{(}\frac{y}{2}\Big{)}\cos\Big{(}\frac{\phi z}{2}\Big{)}+\phi^{-1} x\cos\Big{(}\frac{\phi x}{2}\Big{)}\cos\Big{(}\frac{y}{2\phi}\Big{)}\cos\Big{(}\frac{z}{2}\Big{)}\\
&-&\phi^{-2}y\sin\Big{(}\frac{x}{2}\Big{)}\sin\Big{(}\frac{\phi y}{2}\Big{)}\cos\Big{(}\frac{z}{2\phi}\Big{)}-\phi^{2}y\sin\Big{(}\frac{x}{2\phi}\Big{)}\sin\Big{(}\frac{y}{2}\Big{)}\cos\Big{(}\frac{\phi z}{2}\Big{)}+
\sqrt{5}y\sin\Big{(}\frac{\phi x}{2}\Big{)}\sin\Big{(}\frac{y}{2\phi}\Big{)}\cos\Big{(}\frac{z}{2}\Big{)}\\
&-&\phi^{2}z\sin\Big{(}\frac{x}{2}\Big{)}\cos\Big{(}\frac{\phi y}{2}\Big{)}\sin\Big{(}\frac{z}{2\phi}\Big{)}-
\phi^{-2}z\sin\Big{(}\frac{\phi x}{2}\Big{)}\cos\Big{(}\frac{y}{2\phi}\Big{)}\sin\Big{(}\frac{z}{2}\Big{)}
+\sqrt{5}z\sin\Big{(}\frac{x}{2\phi}\Big{)}\cos\Big{(}\frac{y}{2}\Big{)}\sin\Big{(}\frac{\phi z}{2}\Big{)}.
\end{eqnarray*}  
\normalsize
Then:
\begin{itemize}
\item[1)] the vector field $\mf{I}=(\mf{I}_{x},\mf{I}_{y},\mf{I}_{z})=\mf{V}+\mf{W}$ has an icosahedral symmetry: if we treat $\mf{I}$ as a map $\mathbb{R}^{3}\mapsto\mathbb{R}^{3}$, for any $\zeta\in\mathbb{I}$  one has $\zeta^{-1}\circ\mf{I}\circ\zeta=\mf{I}$; moreover, $\mf{W}$ has the full icosahedral symmetry $\mathbb{I}\times\{I,-I\}$; 
\item[2)] as a Taylor series, $\mf{V}$ contains only terms with even compound degree, $\mf{W}$ contains only odd-degree terms, $\nabla\times\mf{V}=\mf{W}$, $\nabla\times\mf{W}=\mf{V}$; thus $\mf{I}$ satisfies the identity 
\begin{eqnarray*}
\nabla\times\mf{I}=\mf{I};
\end{eqnarray*}
\item[3)] the Taylor series for $\mf{V}$ starts with a degree $6$ vector field $\frac{\mathbf{M}}{768}$, where $\mathbf{M}=(\varpi,\varrho,\sigma)$ is given by \footnotesize
\begin{eqnarray}
\left\{\begin{array}{c@{\qquad}l}
\varpi=(5-\sqrt{5})yz^5+(5+\sqrt{5})y^5z-20y^3z^3+(10+10\sqrt{5})x^2yz^3+
(10-10\sqrt{5})x^2y^3z-10x^4yz,\\
\varrho=(5-\sqrt{5})zx^5+(5+\sqrt{5})z^5x-20z^3x^3+(10+10\sqrt{5})y^2zx^3+(10-10\sqrt{5})y^2z^3x-10y^4zx,\\
\sigma=(5-\sqrt{5})xy^5+(5+\sqrt{5})x^5y-20x^3y^3+(10+10\sqrt{5})z^2xy^3+(10-10\sqrt{5})z^2x^3y-10z^4xy;
\end{array}\right.\label{pirmas}
\end{eqnarray}
\normalsize
 the Taylor series for $\mf{W}$ starts from a degree $5$ vector field $\frac{\m{N}}{768}$, where $\mathbf{N}=\big{(}\lambda,\xi,\chi\big{)}$, $\xi=\lambda(y,z,x)$, $\chi=\lambda(z,x,y)$, and
\footnotesize
\begin{eqnarray*}
\lambda=(35-5\sqrt{5})xy^4-(35+5\sqrt{5})xz^4+60\sqrt{5}xy^2z^2
-(70+10\sqrt{5})x^3y^2+(70-10\sqrt{5})x^3z^2+2\sqrt{5}x^5.
\end{eqnarray*}
\end{itemize}
\label{thm1}
\end{thm}
We can now state one corollary. Suppose, a point $\m{x}\in\mathbb{R}^{3}$, $\m{x}=(x_{0},y_{0},z_{0})\neq\m{0}$, has a non-trivial stabilizer subgroup $G_{\m{x}}<\mathbb{I}$, which is then a rotation of order $2$, $3$, or $5$ with respect to the axis $\mathbb{R}\m{x}$. Then, since a vector $\mf{I}(x_{0},y_{0},z_{0})$ is unchanged under the action of $G_{\m{x}}$, we immediately get that
\begin{eqnarray*}
\mf{I}(x_{0},y_{0},z_{0})=C(x_{0},y_{0},z_{0})\cdot(x_{0},y_{0},z_{0}).
\end{eqnarray*} 
All such points $\m{x}$ with non-trivial stabilizer subgroups lie on $\frac{60}{5}+\frac{60}{3}+\frac{60}{2}=62$ lines, obtained from $\mathbb{R}(\phi,1,0)$ (12 lines corresponding to centres of faces of a dodecahedron), $\mathbb{R}(1,1,1)$ (20 lines corresponding to vertices of a dodecahedron), and $\mathbb{R}(1,0,0)$ (30 lines corresponding to centres of edges) under the action of the group $\mathbb{I}$. In the last section we will see that on every of these lines there exist infinitely many points where the vector field $\mf{I}$ vanishes. For example, if $s_{0}$ is the root of
\begin{eqnarray*}
1-\phi\cos(s)+\phi^{-1}\cos(\phi s)=0,
\end{eqnarray*} 
 then $\mf{I}(\phi s_{0},s_{0},0)=\m{0}$. The smallest positive non-zero such $s_{0}$ is given by $s_{0}=5.1625967944_{+}$. We can call zeros of $\mf{I}$ lying on these $62$ lines as \emph{trivial zeros}.
\begin{prob}
\label{prob-nt}
Are there any non-trivial zeros of the vector field $\mf{I}$?
\end{prob}
Each of these exceptional $62$ lines split into disjoint open segments, each being a complete orbit. The endpoints of the segment are two trivial zeros. The dynamics of all the rest orbits is considerably far more complicated.\\
   
The vector field $\mathbf{M}$ is the numerator of the vector field for the icosahedral projective superflow \cite{alkauskas-super1,alkauskas-super2} (there are two related notions - \emph{projective superflow} and \emph{polynomial superflow}), whence the motivation and the first step comes from.  One has $\zeta^{-1}\circ\m{M}\circ\zeta=\m{M}$ for any $\zeta\in\mathbb{I}$, and this is the unique, up to the scalar multiple, polynomial $6$-homogeneous vector field with this property \cite{alkauskas-super2}. Calculations show that $\nabla\times\m{M}=\m{N}$ has a full icosahedral symmetry, and $\nabla\times(\nabla\times\m{M})=\mathbf{0}$. There exist exactly $5$ irreducible projective superflows in dimension $3$ \cite{alkauskas-super2}. However, the orbits of all five $3$-dimensional superflows - the tetrahedral (group of order $24$, full symmetry), the octahedral ($24$, orientation-preserving symmetry), the icosahedral ($60$, also orientation-preserving symmetry), $3$-prismal ($12$), and $4$-antiprismal ($16$) - are algebraic curves, since the corresponding system of differential equations in all cases possesses two independent algebraic first integrals. This shows a deep analogy (and essential differences) between integration of superflows, and integration of dynamical systems whose Lagrangian has infinitesimal symmetries, exactly as E. Noether's 1918 theorem tells (see \cite{pet}, Chapter IV, \S 12, and also \cite{banados}). We recall that $\mathscr{W}$ is \emph{the first integral} for a vector field $(a_{x},a_{y},a_{z})$, if
\begin{eqnarray*}
\frac{\p\mathscr{W}}{\p x}\cdot a_{x}+\frac{\p\mathscr{W}}{\p y}\cdot a_{y}+\frac{\p \mathscr{W}}{\p z}\cdot a_{z}=0.
\end{eqnarray*}

\begin{figure}
\epsfig{file=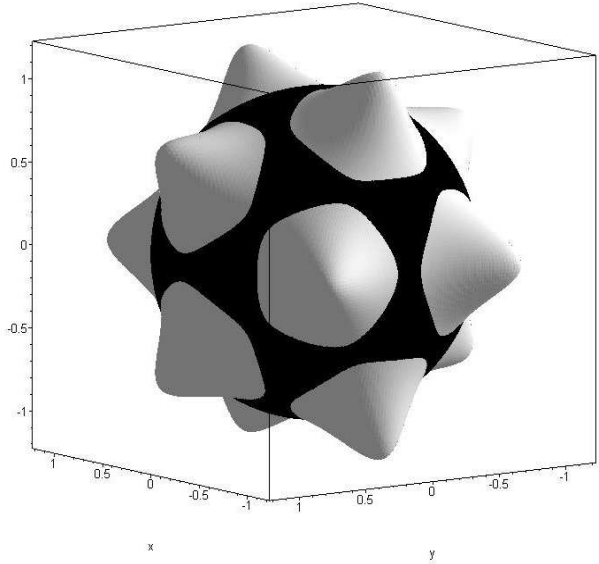,width=270pt,height=250pt,angle=0}
\caption{The intersection of the unit sphere $x^2+y^2+z^2=1$ (black) with the surface  $(\phi^2x^2-y^2)(\phi^2 y^2-z^2)(\phi^2z^2-x^2)+(x^2+y^2+z^2)^3=\frac{19}{20}$ (gray), are separate $12$ orbits of the flow with the vector field $\m{M}$ given by (\ref{pirmas}).}
\label{figure-neg}
\end{figure}

\begin{figure}
\epsfig{file=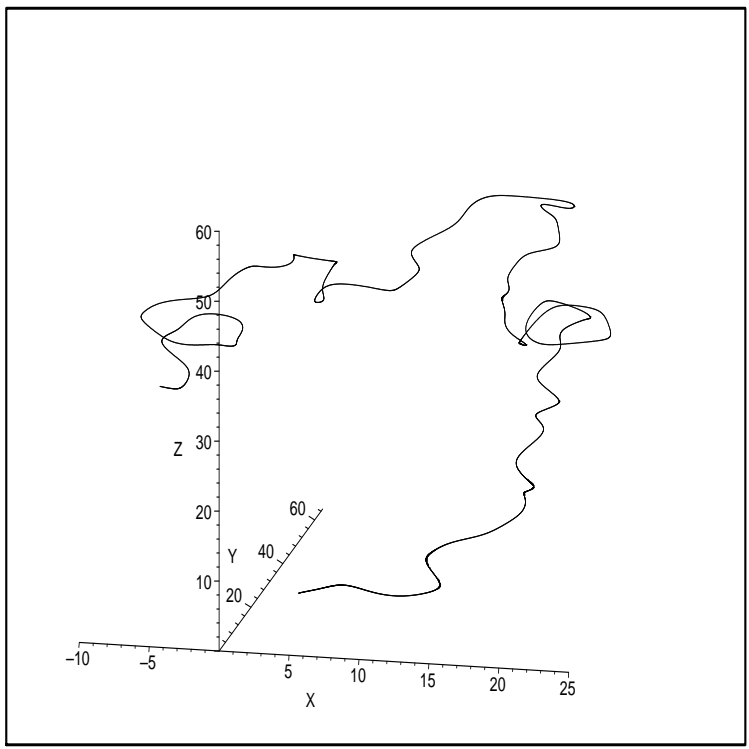,width=280pt,height=280pt,angle=-90}
\caption{Let $\m{x}=(x,y,z)$. Orbit of the point $(5.0,6.0,7.0)$ under the flow $F(\m{x},t)$ with a vector field $\mf{I}$ for the time parameter $0\leq t\leq 1$. In order to numerically calculate the orbit, that is, to solve the system (\ref{system}), we use the classical Runge-Kutta method, which is the fourth order method \cite{iserles}, implemented as a default method in MAPLE to solve a system of ODEs numerically.
}
\label{fig-3}
\end{figure}

Two independent first integrals of the flow with the vector field $\m{M}$ are given by
\begin{eqnarray*}
x^2+y^2+z^2,\quad (\phi^2x^2-y^2)(\phi^2 y^2-z^2)(\phi^2z^2-x^2).
\end{eqnarray*}
 Thus, integration of superflows has a strong algebro-geometric and number-theoretic side \cite{alkauskas-super1, alkauskas-super2,alkauskas-super3}.  For example, $12$ particular orbits of the flow with the vector field $\m{M}$ are shown in Figure \ref{figure-neg}. Note only that differently from the vector field given in \cite{alkauskas-super2}, there is no denominator $(x^2+y^2+z^2)^2$. But since $x^2+y^2+z^2$ is the first integral of the corresponding differential system, that is, $x\varpi+y\varrho+z\sigma\equiv0$, this does not matter. However, in case of the current paper, for the differential system
\begin{eqnarray}
\left\{\begin{array}{c@{\qquad}l}
\dot{x}(t)=\mf{I}_{x}\big{(}x(t),y(t),z(t)\big{)},\\
\dot{y}(t)=\mf{I}_{y}\big{(}x(t),y(t),z(t)\big{)},\\
\dot{z}(t)=\mf{I}_{z}\big{(}x(t),y(t),z(t)\big{)},
\end{array}\right.
\label{system}
\end{eqnarray}
we cannot expect the existence of a single closed-form first integral, not to mention two of them, as Figure \ref{fig-3} suggests. Thus, differently from Noether's theorem for Lagrangians, symmetry is not the only feature of the superflow (with the vector field $\m{M}$, for example) which guarantees the existence of (polynomial) first integrals - superflows is in essence an algebro-geometric topic; for more on this, see \cite{alkauskas-super1}.\\

The terms in the Taylor expansion of $\mf{V}_{x}(x,y,z)$ are not just of even compound degree, but in fact of even degree in $x$ and odd in $y$ and $z$; this is a consequence of the fact $\mathbb{K}<\mathbb{I}$. Cyclically so for $\mf{V}_{x}(y,z,x)$ and $\mf{V}_{x}(z,x,y)$. Equally, $\mf{W}_{x}$ is odd in $x$ and even in each of $y,z$. \\

Note that vector fields $\mf{V}$ and $\mf{W}$ have one additional symmetry. Indeed, let us denote by $\tau$ the non-trivial automorphism of the number field $\mathbb{Q}\big{(}\sqrt{5}\,\big{)}$. Thus, $\tau\phi=-\phi^{-1}$. In the above formulas for $\mf{V}_{x}$ and $\mf{W}_{x}$ (Theorem \ref{thm1}), let us expand everything in Taylor series, swap $y$ and $z$, leave $x$ intact, and apply $\tau$ summand-wise. We readily obtain
\begin{eqnarray*}
\tau\,\mf{V}_{x}(x,z,y)=\mf{V}_{x}(x,y,z), \quad \tau\,\mf{W}_{x}(x,z,y)=-\mf{W}_{x}(x,y,z).
\end{eqnarray*}
The icosahedral group, as a group, is isomorphic to $A_{5}$. The latter group has two non-equivalent $3$-dimensional representations, and the second one is given exactly via an embedding $\tau\alpha=\alpha$, $\tau\beta=\beta$, $\tau\gamma\neq\gamma$ (\cite{kostrikin}, Chapter 8, \S 5, Problem 7).\\

The second result of this paper is completely analogous, only the corresponding Taylor series starts from a degree $9$ rather than $5$, and due to a careful choice of multiplying parameter, there is the same symmetry with respect to $\mathbb{Q}\big{(}\sqrt{5}\,\big{)}$; see (\ref{Y}) further, where any scalar multiple of the collection given produces the solution to our problem, but only $\mathbb{Q}$-scalar multiples have this additional symmetry. The fact that the Taylor coefficients vanish to the order $8$ rather than $4$ makes this vector field even more exceptional. Formulas are more lengthy, so we write down explicitly only the even part.
\begin{thm} 
Let us define the vector field $\mf{V}^{0}=\big{(}\mf{V}^{0}_{x},\mf{V}^{0}_{y},\mf{V}^{0}_{z}\big{)}$, $\mf{V}^{0}_{y}=\mf{V}^{0}_{x}(y,z,x)$, $\mf{V}^{0}_{z}=\mf{V}^{0}_{x}(z,x,y)$, whose Taylor series contains only even compound degrees, by
\begin{scriptsize}
\begin{eqnarray*}
\mf{V}_{x}^{0}&=&2x\sin\Big{(}\frac{x}{2}\Big{)}\sin\Big{(}\frac{\phi y}{2}\Big{)}\sin\Big{(}\frac{z}{2\phi}\Big{)}
-2\phi x\sin\Big{(}\frac{x}{2\phi}\Big{)}\sin\Big{(}\frac{y}{2}\Big{)}\sin\Big{(}\frac{\phi z}{2}\Big{)}+2\phi^{-1}x\sin\Big{(}\frac{\phi x}{2}\Big{)}\sin\Big{(}\frac{y}{2\phi}\Big{)}\sin\Big{(}\frac{z}{2}\Big{)}\\
&+&2\phi y\sin z+(7-\sqrt{5})y\cos\Big{(}\frac{x}{2}\Big{)}\cos\Big{(}\frac{\phi y}{2}\Big{)}\sin\Big{(}\frac{z}{2\phi}\Big{)}\\
&+&2\phi^{2}y\cos\Big{(}\frac{x}{2\phi}\Big{)}\cos\Big{(}\frac{y}{2}\Big{)}\sin\Big{(}\frac{\phi z}{2}\Big{)}
+2\sqrt{5}y\cos\Big{(}\frac{\phi x}{2}\Big{)}\cos\Big{(}\frac{y}{2\phi}\Big{)}\sin\Big{(}\frac{z}{2}\Big{)}\\
&-&2\phi^{-1}z\sin y-(7+\sqrt{5})z\cos\Big{(}\frac{x}{2}\Big{)}\sin\Big{(}\frac{\phi y}{2}\Big{)}\cos\Big{(}\frac{z}{2\phi}\Big{)}\\
&-&2\phi^{-2}z\cos\Big{(}\frac{\phi x}{2}\Big{)}\sin\Big{(}\frac{y}{2\phi}\Big{)}\cos\Big{(}\frac{z}{2}\Big{)}
-2\sqrt{5}z\cos\Big{(}\frac{x}{2\phi}\Big{)}\sin\Big{(}\frac{y}{2}\Big{)}\cos\Big{(}\frac{\phi z}{2}\Big{)},
\end{eqnarray*}
\end{scriptsize}
and its odd counterpart $\mf{W}^{0}$ is defined by $\mf{W}^{0}=\nabla\times\mf{V}^{0}$. 
Then the vector field $\mf{Y}=\mf{V}^{0}+\mf{W}^{0}$ has an icosahedral symmetry, and satisfies $\nabla\times\mf{Y}=\mf{Y}$. The Taylor series for $\mf{V}^{0}$ starts from the degree $10$ vector field $\frac{\m{P}}{23224320}$, where $\m{P}=(\varpi_{0},\varrho_{0},\sigma_{0})$, $\varrho_{0}(x,y,z)=\varpi_{0}(y,z,x)$, $\sigma_{0}(x,y,z)=\varpi_{0}(z,x,y)$, and  where\footnotesize
\begin{eqnarray}
\varpi_{0}&=&-18x^{8}yz+(84+84\sqrt{5})x^{6}y^{3}z+(84-84\sqrt{5})
x^{6}yz^{3}-(126+126\sqrt{5})x^{4}y^{5}z-(126-126\sqrt{5})x^{4}yz^{5}
\nonumber\\
&+&(36+108\sqrt{5})x^{2}y^{7}z+(36-108\sqrt{5})x^2yz^{7}-504\sqrt{5}x^{2}y^{5}z^{3}
+504\sqrt{5}
x^{2}y^{3}z^{5}
\label{antras}\\
&+&(9-5\sqrt{5})y^{9}z+(9+5\sqrt{5})yz^{9}-(120-24\sqrt{5})y^{7}z^{3}
-(120+24\sqrt{5})y^{3}z^{7}+252y^{5}z^{5}.
\nonumber
\end{eqnarray}
\normalsize
The Taylor series for the vector field $\mf{Y}$ starts from a degree $9$ vector field 
\begin{eqnarray*}
\frac{\m{Q}}{23224320}=\nabla\times\frac{\m{P}}{23224320}\neq\m{0}.
\end{eqnarray*}
\label{thm2}
\end{thm} 
Thus, a one-parameter family of vector fields $\mf{I}_{a}=\mf{I}+a\mf{Y}$ has $\frac{\mathbf{N+M}}{768}$ as a beginning of its Taylor series expansion, has an icosahedral symmetry, and satisfies $\nabla\times\mf{I}_{a}=\mf{I}_{a}$. If $a\in\mathbb{Q}$, the same symmetry with respect to $\mathbb{Q}(\,\sqrt{5}\,)$ applies. As we will soon see in Note in Section \ref{icos}, this family arises from a $1$-dimensional setting for the \emph{Helmholtz equation}. This is what we meant by saying \emph{the simplest possible} in the very introduction to this paper. Higher dimensional solutions to the Helmholtz equation, which lead to Beltrami vector fields with various polyhedral symmetries, are treated in \cite{alkauskas-beltrami}; see Section \ref{order}.\\

If two $n$-dimensional vector fields $X$ and $Y$ in $\mathbb{R}^{n}$ are given by
\begin{eqnarray*}
X=\sum\limits_{i=1}^{n}f_{i}\frac{\p}{\p x_{i}},\quad 
Y=\sum\limits_{i=1}^{n}g_{i}\frac{\p}{\p x_{i}},
\end{eqnarray*}
then their \emph{Lie bracket} is defined as \cite{conlon,gadea,hitchin}
\begin{eqnarray*}
[X,Y]=\sum\limits_{i=1}^{n}\Big{(}X(g_{i})-Y(f_{i})\Big{)}\frac{\p}{\p x_{i}}.
\end{eqnarray*}
The Taylor series of $[\mf{I},\mf{Y}]$ starts from a scalar multiple of a vector field $[\m{N},\m{Q}]$. Computer calculations show that this is not identically $0$, and so
\begin{corr}Vector fields $\mf{I}$ and $\mf{Y}$ do not commute.
\end{corr}
\section{Overview}
\label{overview}
One of the motivations to investigate Beltrami fields is the following  claim which follows from the result by Arnold \cite{arnold,etnyre}, and which shows that Beltrami flows on a closed $3$-manifold might have a complicated topology. Namely, the orbits of a Beltrami flow on a $3$-manifold are not always constrained to lie on a $2$-torus as is the case for steady Euler flows with a $C^{\omega}$ vector field and which are not everywhere collinear with its curl. For more information on various aspects of Beltrami flows, see \cite{blair,hiptmair,macleod1,macleod2,torres,woltjer}. For mathematical physics-related aspects of force-free magnetic fields, see \cite{chandr1,chandr2,chandr3,kaiser, kravchenko,priest,salingaros,wiegelmann}.\\

It is impossible to give a wide overview on a prolific literature related to Beltrami fields, so we will confine to few papers. \\

For a computational aspect (in bounded simply-connected domains), see \cite{amari}. Symmetry questions of force-free fields not depending on the variable $z$ (the so called $2$-dimensional FFF) are treated in \cite{tassi}. Recall that \emph{a planefield} $\xi$ on a $3$-dimensional manifold is a smooth mapping that assigns to every point $p$ a plane in its tangent space. A planefield $\xi$ is said to be \emph{integrable} at a point $p$, if there exists a smooth surface $S$ passing through $p$ such that $\xi$ is tangential to $S$ in some neighbourhood of $p$. A planefield $\xi$ is called \emph{a contact structure} if and only if it is everywhere non-integrable. The relation between Beltrami (and Trkalian) fields with contact structures is investigated in \cite{dahl,kom1,kom2}. A topic, slightly related to our paper, was investigated in \cite{brandolese}, where it was shown that non-stationary solution of the Navier-Stokes equations in $\mathbb{R}^{d}$ ($d=2,3$), if this solution is left invariant under the action of a finite subgroup of the orthogonal group, decays much faster as $|x|\rightarrow\infty$ or $t\rightarrow\infty$ than in a generic case. The decay is extremely fast in case $d=3$ and the full symmetry group of an icosehedron; that is, $\mathbb{I}\times\{I,-I\}$. In \cite{kaiser} the authors prove, via an iteration scheme, the existence of force-free magnetic fields in the exterior domain of some compact simply connected surface $S$. The huge difference emerges if we consider constant, or non-constant force-free fields. The spherical curl transform for Trkalian fields using differential forms and its Radon transform are investigated in \cite{saygili1,saygili2}. In \cite{morse1,morse2} the author investigates the eigenfunctions of the equation $\nabla\times\m{B}=\lambda\m{B}$ for finite cylindrical geometry with normal boundary condition $\vec{n}\cdot\m{B}=0$ for nonaxisymmetric modes. The author investigates also the equation $\nabla\times\nabla\times\m{B}=\lambda^{2}\m{B}$, since this reveals the underlying elliptic nature of the initial eigenvector problem. The method of Chandrasekhar and Kendall (which we will also employ) is being used. Double-curl spectral Beltrami equation
\begin{eqnarray*}
\left\{\begin{array}{c@{\qquad}l}
\nabla\times(\nabla\times\m{B})+
\alpha\nabla\times\m{B}+\beta\m{B}=0 \text{ (in the domain }\Omega),\\
\vec{n}\cdot\m{B}=0,\quad \vec{n}\cdot(\nabla\times\m{B})=0\text{ (on }\partial\Omega),
\end{array}\right.
\end{eqnarray*}
are investigated in \cite{mahajan1, mahajan2}. The authors show that if the domain $\Omega$ is multiply connected, then this equation has a nonzero solution for arbitrary
complex  numbers $\alpha$ and $\beta$.\\

In relation to the results of the current paper, we will emphasize one consequence of the result proved in \cite{yoshida} (see \cite{comp} for numerical results in this direction).\\

Consider the simplest example of a Beltrami condition satisfied by a three-dimensional solenoidal
vector field, obeying
\begin{eqnarray}
\left\{\begin{array}{c@{\qquad}l}
\nabla\times\m{B}=\lambda\m{B}\text{ (in the domain }\Omega),\\
\vec{n}\cdot\m{B}=0\text{ (on }\partial\Omega).
\end{array}\right.
\label{eigen}
\end{eqnarray}
Here $\lambda$ is a real (or complex) constant number, $\Omega\subset\mathbb{R}^{3}$ is a bounded domain with a smooth boundary, and $\vec{n}$ is a unit normal vector onto $\Omega$. This system is regarded as an eigenvalue problem with respect to the curl operator. Then one of the results in \cite{yoshida} claims that if $\Omega$ is simply connected, then the system (\ref{eigen}) has a nonzero solution for special $\lambda$ included in a set of discrete real numbers; these numbers represent the point spectrum of the self-adjoint part of the curl operator.\\

Now, consider a simply-connected domain $\Omega$ which has an icosahedral symmetry $\mathbb{I}$. For example, we can take a domain (see Figure \ref{figure-neg})
\begin{eqnarray}
\Omega=\{(x,y,z)\in\mathbb{R}^{3}:(\phi^2x^2-y^2)(\phi^2 y^2-z^2)(\phi^2z^2-x^2)+(x^2+y^2+z^2)^3
< 1\}.
\label{star}
\end{eqnarray}
Let $(\m{B},\lambda)$ solves an eigenvalue problem (\ref{eigen}). It is obvious that any orientation-preserving  orthogonal change of coordinates $\gamma^{-1}\circ\m{B}\circ\gamma(x,y,z)$, $\gamma\in\mathbb{I}$, also solves this problem. Thus,
\begin{eqnarray}
\m{B}_{\mathbb{I}}=\frac{1}{60}\sum\limits_{\gamma\in\mathbb{I}}\gamma^{-1}\circ\m{B}\circ{\gamma}
\label{average}
\end{eqnarray}
is a vector field with an icosahedral symmetry which also solves the eigenvalue problem (\ref{eigen}). However, without delving deeper into the geometry of the surface $\partial\Omega$ and the system (\ref{eigen}), we cannot guarantee that $\m{B}_{\mathbb{I}}$ is non-zero. Further, as was shown in the introduction, there exist $62$ lines (see also Section \ref{paskut}) passing through the origin such that for points on these lines, a vector field $\m{B}_{\mathbb{I}}$ is collinear with a line itself. Suppose now that a domain $\Omega$ is a star domain with respect to the origin, which implies that intersection of every line through the origin with $\overline{\Omega}$ is a closed interval. The domain (\ref{star}) is an example. Then we have the following result. (By icosahedral symmetry we always mean the group $\mathbb{I}$).
\begin{prop}If a vector field $\m{B}_{\mathbb{I}}$ has an icosahedral symmetry, and $(\m{B}_{\mathbb{I}},\lambda)$ solves the eigenvalue problem (\ref{eigen}) for the star domain (with respect to the origin) $\Omega$ with an icosahedral symmetry, then there exist at least $62$ points on $\Omega$ where a vector field $\m{B}_{\mathbb{I}}$ vanishes.
\end{prop}

This gives a new perspective on uniqueness of vector fields $\mf{I}$ and $\mf{Y}$ which are constructed in this paper.\\
  
Next, we formulate one problem which seems to be of a huge interest from the point of view of geometry. In \cite{blair}, Section 7.2, the author shows that the existence of conformally flat contact metric manifolds corresponds to finding the solutions to the equation (in cartesian coordinates) $\nabla\times\m{B}=|\m{B}|\m{B}$, where $|\m{B}|$ is vectors length. For the standard Sasakian structure of a constant curvature $+1$ on $\m{S}^{3}$, using stereographic projection to $\mathbb{R}^{3}$, the corresponding vector field is 
\begin{eqnarray}
\m{B}=\frac{8(xz-y)}{(1+x^2+y^2+z^2)^2}\frac{\p}{\p x}
+\frac{8(x+yz)}{(1+x^2+y^2+z^2)^2}\frac{\p}{\p y}
+\frac{4(1+z^2-x^2-y^2)}{(1+x^2+y^2+z^2)^2}\frac{\p}{\p z}.
\label{vec-b}
\end{eqnarray}
Here 
\begin{eqnarray}
|\m{B}|=\frac{4}{1+x^2+y^2+z^2}.
\end{eqnarray} 
In relation to this, as a more algebro and differential-geometric version of results of the current paper, we pose
\begin{prob}
\label{prob-rb}
Does there exist a Beltrami vector field $\m{B}$, meaning $\nabla\times\m{B}=f(x,y,z)\cdot\m{B}$, which is given by rational functions, and which has a 4-Klein group? tetrahedral? octahedral? icosahedral symmetry?
\end{prob}
Note that, if $\nabla\times\m{B}=f(x^2+y^2+z^2)\cdot\m{B}$, then
\begin{eqnarray*}
\nabla\times(\eta^{-1}\circ\m{B}\circ\eta)=f(x^2+y^2+z^2)\cdot\eta^{-1}\circ\m{B}\circ\eta
\end{eqnarray*} 
for any $\eta\in SO(3)$. Therefore, it as if seems that we can construct Beltrami fields with any cyclic or polyhedral symmetry from any given Beltrami field by averaging, like (\ref{average}). However, in most cases we will end up with a $0$ vector field, and independent methods are truly needed - this is the whole essence of this paper! For example, the Klein four group
\begin{eqnarray}
\mathbb{K}=\{I,\mathrm{diag}(-1,-1,1),\mathrm{diag}(-1,1,-1),\mathrm{diag}(1,-1,-1)\}
\label{klein}
\end{eqnarray}
is a subgroup of $\mathbb{I}$. But yet, for $\m{B}$ given by (\ref{vec-b}),
\begin{eqnarray*}
\sum\limits_{\eta\in\mathbb{K}}\eta^{-1}\circ\m{B}\circ\eta=\m{0}.
\end{eqnarray*}
Still, if $\beta$ is given by (\ref{gen-ico}), we get a non-trivial example by calculating
\begin{eqnarray*}
\m{F}=\frac{1}{4}\sum\limits_{j=0}^{2}\beta^{-j}\circ\m{B}\circ\beta^{j}=\big{(}U(x,y,z),U(y,z,x),U(z,x,y)\big{)}
\end{eqnarray*} 
(factor $4$, not $3$, is for simplicity), where
\begin{eqnarray*}
U=\frac{2xy+2xz-2y+2z+1+x^2-y^2-z^2}{(1+x^2+y^2+z^2)^2}.
\end{eqnarray*}
Thus, the answer to Problem \ref{prob-rb} is positive at least in a cyclic order $3$ subgroup generated by $\beta$. Moreover, we have
\begin{eqnarray*}
\nabla\times\m{F}=\frac{4}{1+x^2+y^2+z^2}\cdot\m{F}=\frac{4}{\sqrt{3}}\cdot|\m{F}|\cdot\m{F}.
\end{eqnarray*}
We will address this problem in the next publication.\\
\section{Construction}
\label{icos}
One of the main identities of a $3$-dimensional vector calculus claims that for a smooth vector field $\m{B}$,
\begin{eqnarray} 
\nabla\times(\nabla\times\m{B})=\nabla(\nabla\cdot\m{B})-\nabla^{2}\m{B};
\label{vec-bas}
\end{eqnarray} 
here $\nabla^{2}$ is a vector Laplace operator, and for a scalar function $f$, $\nabla f=\mathrm{grad}\, f$. To prove our two theorems, first we will construct a vector field 
$\mf{V}$ such that:
\begin{itemize}
\item[1)] $\mf{V}$ satisfies the vector Helmholtz equation $\nabla^{2}\mf{V}=-\mf{V}$; 
\item[2)] it satisfies $\nabla\cdot\mf{V}=0$; 
\item[3)] $\mf{V}$ has an icosahedral symmetry;
\item[4)] all elements in the Taylor series are of even compound degree.
\end{itemize} 
If the first two are satisfied, then the identity (\ref{vec-bas}) implies
\begin{eqnarray*}
\nabla\times(\nabla\times\mf{V})=\mf{V}.
\end{eqnarray*} 
Then we put $\mf{W}=\nabla\times\mf{V}$, we will see that $\mf{I}=\mf{V}+\mf{W}$ has the properties described by items 1) and 2) in Theorem \ref{thm1}. Such method of deriving solutions to $\nabla\times\m{B}=\nu\m{B}$ from \emph{scalar} solutions to the Helmholtz equation $(\nabla^{2}+\nu^{2})\Psi=0$ (then such solution $\Psi$ is called \emph{Debye potential}) was developed by Chandresekhar and Kendall \cite{chandr1,saygili1}. Our contribution is the fact that we consider vector solutions, and especially the emphasis on part 3), which is new in the theory of Beltrami fields. See the end of Section \ref{sec-comp} where it is shown that orthogonal orientation-preserving symmetries of a vector field carries automatically as symmetries of its curl. \\

To satisfy the requirements 1), 2) and 3) above (we now secure the notation $\mf{V}$ for a specific vector field given by Theorem \ref{thm1}, changing the unspecified vector field to $\mf{H}$), we will try to find constants $a_{i}$, and homogeneous linear forms $L_{i}$, $k_{i}$, such that $\mf{H}=(\mf{G}_{x},\mf{G}_{y},\mf{G}_{z})$, where
\begin{eqnarray}
\mf{G}=\mf{G}_{x}=\sum\limits_{i}\big{(}a_{i}\cos(k_{i} )+L_{i}\sin(k_{i})\big{)},
\label{ico-sum}
\end{eqnarray} 
and the two other coordinates $\mf{G}_{y}$ and $\mf{G}_{z}$ are obtained from $\mf{G}$ by a cyclic permutation. First, we know that $\mf{H}=(\mf{G}_{x},\mf{G}_{y},\mf{G}_{z})$ is invariant under conjugation with $\alpha=\mathrm{diag}(-1,-1,1)$, $\beta\alpha\beta^{2}=\mathrm{diag}(1,-1,-1)$, and $\beta^{2}\alpha\beta=\mathrm{diag}(-1,1,-1)$. These four matrices, as already mentioned two times, together with the unity $I$ produce the Klein's fourth group $\mathbb{K}<\mathbb{T}<\mathbb{I}$. This gives 
\begin{eqnarray}
-\mf{G}(-x,-y,z)=\mf{G}(x,y,z),\quad\mf{G}(x,-y,-z)=\mf{G}(x,y,z).
\label{inv-k}
\end{eqnarray} 
As the first three linear forms $k_{i}$, we just take $x,y,z$. Invariance under $\mathbb{K}$ gives, respectively, the following sum as a part of $\mf{G}_{x}$ in (\ref{ico-sum}); namely,
$az\sin y+by\sin z$. Next, let us define the following $12$ linear forms and vectors as follows.
\begin{eqnarray*}
\begin{tabular}{l l}
$\ell_{x,0}(x,y,z)=\frac{1}{2}x+\frac{\phi}{2}y+\frac{1}{2\phi}z$,& $\m{j}_{x,0}=\Big{(}\frac{1}{2},\frac{\phi}{2},\frac{1}{2\phi}\Big{)},$\\
$\ell_{x,1}(x,y,z)=-\frac{1}{2}x+\frac{\phi}{2}y+\frac{1}{2\phi}z$,&
$\m{j}_{x,1}=\Big{(}-\frac{1}{2},\frac{\phi}{2},\frac{1}{2\phi}\Big{)},$\\
$\ell_{x,2}(x,y,z)=\frac{1}{2}x-\frac{\phi}{2}y+\frac{1}{2\phi}z$,&
$\m{j}_{x,2}=\Big{(}\frac{1}{2},-\frac{\phi}{2},\frac{1}{2\phi}\Big{)},$\\
$\ell_{x,3}(x,y,z)=\frac{1}{2}x+\frac{\phi}{2}y-\frac{1}{2\phi}z$,&
$\m{j}_{x,3}=\Big{(}\frac{1}{2},\frac{\phi}{2},-\frac{1}{2\phi}\Big{)}$.
\end{tabular}
\end{eqnarray*}

Similarly we define $8$ other linear functions and $8$ vectors by cyclically permuting variables. For example,  $\m{j}_{y,2}=(\frac{1}{2\phi},\frac{1}{2},-\frac{\phi}{2})$, $\ell_{y,2}=\m{j}_{y,2}\cdot(x,y,z)^{T}=\frac{1}{2}y-\frac{\phi}{2}z+\frac{1}{2\phi}x$, and so on. Orthogonality of the matrix $\gamma$ tells, for example, that $\m{j}_{x,2}$, $\m{j}_{z,1}$ and $\m{j}_{y,0}$ are orthonormal vectors. All these relations amount to the same identities $\phi^{2}+\phi^{-2}+1=4$ (unit length), or $\phi-\phi^{-1}-1=0$ (orthogonality). This gives $14$ linear forms in total ($12$ plus two forms $y$ and $z$), and this is our complete collection in (\ref{ico-sum}).\\

Next, if the coordinate of the form (\ref{ico-sum}) is invariant under conjugating with $\mathbb{K}$, it is of the form  \small
\begin{eqnarray}
\mf{G}&=&az\sin y+by\sin z
+c\cos\ell_{x,0}-c\cos\ell_{x,3}-c\cos\ell_{x,2}+c\cos\ell_{x,1}\label{one-lost}\\
&+&K(x,y,z)\sin\ell_{x,0}+K(-x,-y,z)\sin\ell_{x,3}+K(-x,y,-
z)\sin\ell_{x,2}-K(x,-y,-z)\sin\ell_{x,1}\nonumber\\
&+&d\cos\ell_{y,0}-d\cos\ell_{y,2}-d\cos\ell_{y,1}+d\cos\ell_{y,3}\nonumber\\
&+&L(x,y,z)\sin\ell_{y,0}+L(-x,-y,z)\sin\ell_{y,2}+L(-x,y,-z)\sin\ell_{y,1}-L(x,-y,-z)\sin\ell_{y,3}\nonumber\\
&+&e\cos\ell_{z,0}-e\cos\ell_{z,1}-e\cos\ell_{z,3}+e\cos\ell_{z,2}\nonumber\\
&+&M(x,y,z)\sin\ell_{z,0}+M(-x,-y,z)\sin\ell_{z,1}+M(-x,y,-z)\sin\ell_{z,3}-M(x,-y,-z)\sin\ell_{z,2}.\nonumber
\end{eqnarray}\normalsize
Here $a,b,c,d,e$ are arbitrary constants, and $K,L,M$ are arbitrary linear forms. \\

Now, recall that we want $\nabla^{2}\mf{G}=-\mf{G}$ to be satisfied. We will achieve this if each of the summands in (\ref{one-lost}) satisfies the Helmholtz equation. Now, the following Lemma is immediate.
\begin{lem}
\label{lemu}
Let $\m{a},\m{b}\neq\m{0}$ be two $3$-dimensional vectors-rows, and $\m{x}=(x,y,z)^{T}$. The function $\m{a}\m{x}\cdot\sin(\m{b}\m{x})$ is a solution to the Helmholtz equation if and only if $\langle\m{a},\m{b}\rangle=0$, and $|\m{b}|=1$. The same holds for the $\cos$ function.
\end{lem}
By a direct inspection, two vectors orthogonal to $\m{j}_{x,0}$ are given by $\m{j}_{z,2}$ and $\m{j}_{y,1}$. So,
\begin{eqnarray*}
K(x,y,z)=f\ell_{z,2}+g\ell_{y,1},\quad
L(x,y,z)=h\ell_{x,2}+i\ell_{z,1},\quad
M(x,y,z)=j\ell_{y,2}+k\ell_{x,1}.
\end{eqnarray*}
We therefore have $11$ free coefficients ($a$ through $k$) at our disposition. Such a function $\mf{G}(x,y,z)$ is invariant under conjugation with $\mathbb{K}$ and satisfies the Helmholtz equation. In vector terms, the vector field $\big{(}\mf{G}(x,y,z),\mf{G}(y,z,x),\mf{G}(z,x,y)\big{)}$ has a tetrahedral symmetry $\mathbb{T}$ of order $12$ (generated by matrices $\alpha$ and $\beta$) and satisfies the vector Helmholtz equation.  We will reduce the amount of free coefficients by requiring that a vector field has also a $\gamma$-symmetry, and that the divergence vanishes. 
\section{Order of approach}
\label{order}
The expression (\ref{one-lost}) is what we mean by the \emph{first order approach}. The second order functions $(x^2-y^2)\cos z$ and $2xy\cos z$, for example, also satisfy the Helmholtz equation. In general, let $n\in\mathbb{N}_{0}$, and $P_{n}(x,y)=\Re(x+iy)^{n}$, $Q_{n}(x,y)=\Im(x+iy)^{n}$ be the standard harmonic polynomials of order $n$. Then the $n$th order solutions to the Helmholtz equation are given by $P_{n}(x,y)\cos z$, $Q_{n}(x,y)\cos z$. Other solutions are given by $\sin z$ instead of $\cos z$, or from any of these by an orthogonal change of variables, and all possible linear combinations. This case is treated in \cite{alkauskas-beltrami}, including a construction of Beltrami  vector fields with tetrahedral and octahedral symmetries.\\

Now, fix $N\in\mathbb{N}$, and consider the expression (\ref{ico-sum}), where this time:
\begin{itemize}
\item[i)]Linear forms $k_{i}$ are given by $12$ linear forms $\ell_{w,a}$, $w\in\{x,y,z\}$, $a\in\{0,1,2,3\}$, and 3 linear forms $x,y,z$;
\item[ii)]$a_{i}$ is a polynomial in $(x,y,z)$ of even compound degree and $L_{i}$ is of odd, respectively, both are of degrees at most $N$;
\item[iii)]Each term in (\ref{ico-sum}) satisfies the Helmholtz equation;
\item[iv)] $(\mf{G}(x,y,z),\mf{G}(y,z,x),\mf{G}(z,x,y))$
has the icosahedral symmetry $\mathbb{I}$;
\item[v)] $(\mf{G}(x,y,z),\mf{G}(y,z,x),\mf{G}(z,x,y))$ is solenoidal.
\end{itemize}
\begin{prob}
\label{prob-is}
Find the dimension $d_{\mathbb{I}}(N)$ of linear space of all such functions $\mf{G}$.
\end{prob}
In the next section we will show tat $d_{\mathbb{I}}(0)=0$, and $d_{\mathbb{I}}(1)=2$.\\

While dealing with the icosahedral group, we are working with the dimension $n=3$. The same set of ideas carries to any dimension, with a slight attenuation of requirements - there is no same dimensional analogue of a curl operator in dimensions other than $3$ ($n=\frac{n(n-1)}{2}$ holds for $n=3$ only, $n\in\mathbb{N}$).\\

 Namely, the property $\nabla\times\mf{G}=\mf{G}$ is replaced by the properties 1), 2) and 4) in the beginning of Section \ref{icos}. This leads, for example, to the vector field $\mathfrak{D}=(\mf{T}_{x},\mf{T}_{y})$, given by
\begin{eqnarray*}
\mf{T}_{x}&=&-\cos y+\sqrt{3}\sin\Big{(}\frac{x}{2}\Big{)}\sin\Big{(}\frac{\sqrt{3}y}{2}\Big{)}+\cos\Big{(}\frac{\sqrt{3}x}{2}\Big{)}\cos\Big{(}\frac{y}{2}\Big{)},\\
\mf{T}_{y}&=&-\cos x+\sqrt{3}\sin\Big{(}\frac{y}{2}\Big{)}\sin\Big{(}\frac{\sqrt{3}x}{2}\Big{)}+\cos\Big{(}\frac{\sqrt{3}y}{2}\Big{)}\cos\Big{(}\frac{x}{2}\Big{)},
\end{eqnarray*}
with these properties.
\begin{itemize}
\item[i)]The vector field $\mathfrak{D}$ has a $6$-fold dihedral symmetry, generated by the  matrices
\begin{eqnarray*} 
\begin{pmatrix}0 & 1\\1 & 0\\ \end{pmatrix},\quad
\frac{1}{2}\begin{pmatrix}-1 & -\sqrt{3}\\
\sqrt{3}& -1\\ \end{pmatrix}.
\end{eqnarray*}
\item[ii)]the Taylor series for $\mathfrak{D}$ contains only even compound degrees, and it starts from
\begin{eqnarray*}
\frac{3}{8}\Big{(}(2xy-x^2+y^2),(2xy+x^2-y^2)\Big{)};
\end{eqnarray*}
\item[iii)]it satisfies the vector Helmholtz equation $\nabla^{2}\mathfrak{D}=-\mathfrak{D}$;
\item[iv)]$\mathrm{div}\,\mathfrak{D}=0$.
\end{itemize}
So, these properties relate the vector field to $ABC$ flows. The current paper is thus an introduction to a much broader setting dealt with in \cite{alkauskas-beltrami}.\\
 
 We remind that \emph{the ABC flow}, or \emph{Arnold-Beltrami-Childress flow}, is described by a vector field 
\begin{eqnarray*}
\mathbf{b}=\Big{(}A\sin z+C\cos y,B\sin x+A\cos z,C\sin y+B\cos x\Big{)},
\end{eqnarray*} 
$A,B,C\in\mathbb{R}$. Symmetry related question of these are treated in \cite{ABC}. Such a vector fields satisfies $\mathbf{b}=\nabla\times\mathbf{b}$, and from the Helmholtz equation point of view, it represents an order $0$ approach. Generally, it does not have any non-trivial orthogonal symmetries, except in some cases, like $A=B=C=1$, but, of course, there many more profound symmetries apart from those generated by the group $2\pi\cdot\mathbb{Z}^{3}$ \cite{ABC}.\\

\section{Computer-assisted calculation}
\label{sec-comp}
Now, we will proceed with determining $11$ free coefficients in (\ref{one-lost}). \\

 \begin{figure}[h]
\centering
\epsfig{file=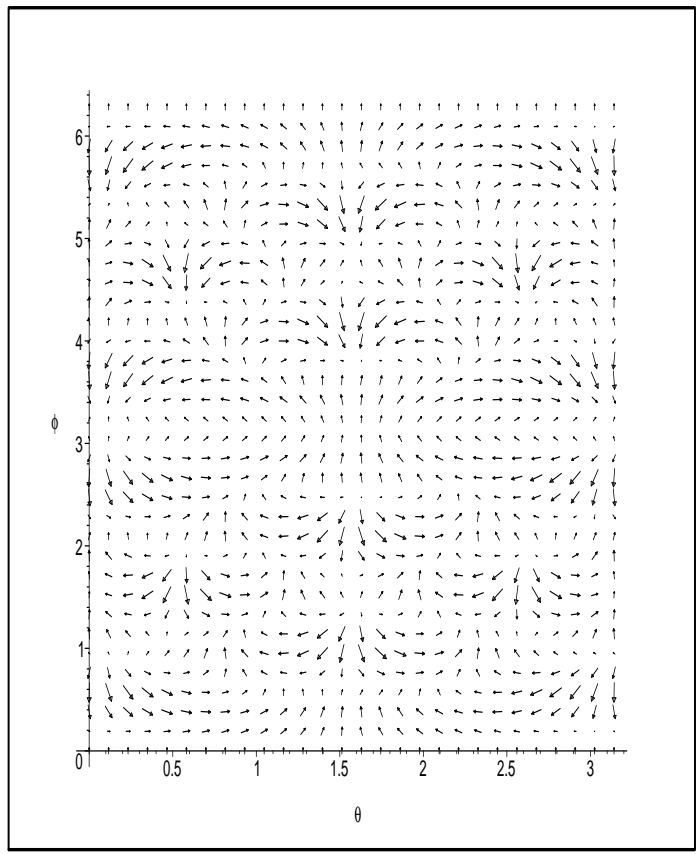,width=320pt,height=260pt,angle=-90}
\caption{Vector field $\mf{I}$ is converted to spherical coordinates $(r,\varphi,\theta)$ as $\big{(}A(r,\varphi,\theta),B(r,\varphi,\theta),C(r,\varphi,\theta)\big{)}$. The plots show vector fields $\big{(}C(1,\varphi,\theta),A(1,\varphi,\theta)\big{)}$ in the $``\theta-\varphi"$ plane, for $(\theta,\varphi)\in[0,\pi]\times[0,2\pi]$ (the unit sphere).}
\end{figure}
 \begin{figure}[h]
\centering
\epsfig{file=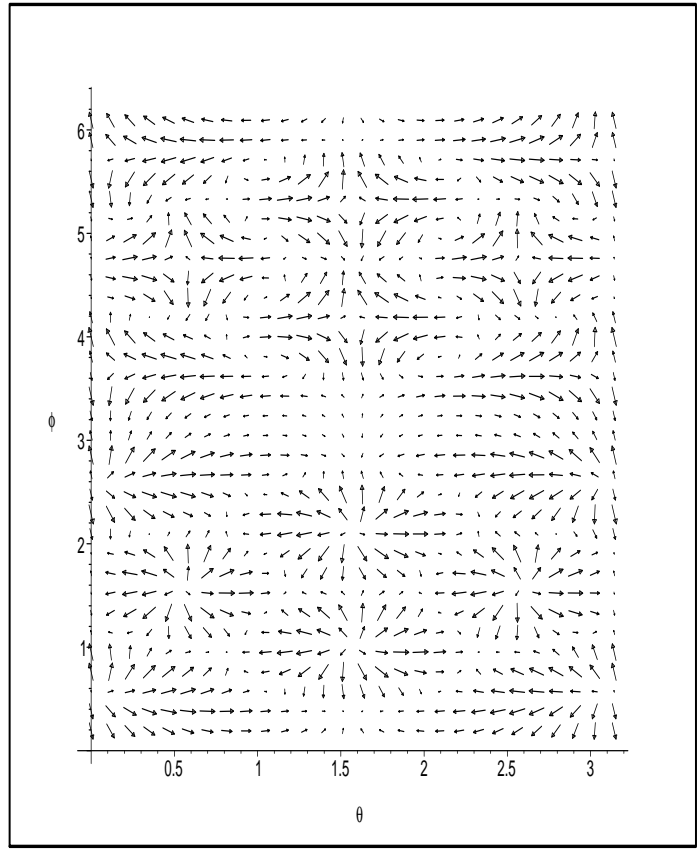,width=320pt,height=260pt,angle=-90}

\caption{The same, only the vector field $\big{(}C(1,\varphi,\theta),B(1,\varphi,\theta)\big{)}$.}
\end{figure}

\indent\textbf{i) }First, we will require that the Taylor coefficient of $\mf{G}$ of degree $6$ is exactly equal to $\varpi$ as given by (\ref{pirmas}), and Taylor coefficients of degrees $\leq 5$ are absent. This gives, with computer calculations performed on MAPLE, $6$ linear relations among $11$ constants. The free coefficients turn out to be $a,b,c,d$, and $f$. This is, of course, to some extent a loose choice; one can also arrive at $5$ other parameters, other $6$ being linear expressions in former ones. This way we obtain a $5$-parameter function $\mf{G}$ given by (\ref{one-lost}), such that its Taylor series starts with $\varpi$, it is invariant under conjugation with $\mathbb{K}$, that is, has a symmetry (\ref{inv-k}), and satisfies the Helmholtz equation. \\

\indent\textbf{ii) }As the next step, we require that $\nabla\cdot\big{(}\mf{G}(x,y,z),\mf{G}(y,z,x),\mf{G}(z,x,y)\big{)}=0$. This gives two more relations, leaving $a,b$, and $d$ as free coefficients. In fact, these computations, though they involve only manipulations with polynomials and Taylor coefficients of trigonometric functions, take some time on a modern computer. It is infeasible to do it by hand. (In \cite{alkauskas-beltrami}, however, we will show that there is an alternative method to calculate the vector field in Theorem \ref{thm1} manually).\\

\indent\textbf{iii) }Finally, a $3$-parameter vector field $\mf{H}=\big{(}\mf{G}(x,y,z),\mf{G}(y,z,x),\mf{G}(z,x,y)\big{)}$ has the tetrahedral symmetry, satisfies the Helmholtz equation, has a vanishing divergence, and the correct beginning of the Taylor series. If it is invariant under conjugation with $\gamma$, it has the icosahedral symmetry, and the problem is solved. We have:
\begin{eqnarray*}
\mf{H}&=&\gamma^{-1}\circ\Big{(}\mf{G}(x,y,z),\mf{G}(y,z,x),\mf{G}(z,x,y)\Big{)}\circ\gamma(x,y,z)\\
&=&\gamma^{-1}\circ\Big{(}\mf{G}(\ell_{x,2},\ell_{z,1},\ell_{y,0}),\mf{G}(\ell_{z,1},\ell_{y,0},\ell_{x,2}),\mf{G}(\ell_{y,0},\ell_{x,2},\ell_{z,1})\Big{)}:=\gamma^{-1}(A,B,C).
\end{eqnarray*}
Since $\m{j}_{x,0}$ is the first row of $\gamma^{-1}=\gamma^{T}$, therefore, we finally require that
\begin{eqnarray*}
\frac{A}{2}+\frac{\phi B}{2}+\frac{C}{2\phi}=\mf{G}(x,y,z).
\end{eqnarray*}

This gives two more relations, and we thus obtain a $1$-parameter family of icosahedral vector fields, with a free parameter being $a$. We will get a particularly symmetric, with respect to a non-trivial automorphism of $\mathbb{Q}\big{(}\sqrt{5}\,\big{)}$, example if $a=\tau b$. Indeed, the last of the $10$ linear equations we obtained reads as $b=1920-\frac{3}{2}a-\frac{\sqrt{5}}{2}a+384\sqrt{5}$. Let us choose $a=384\cdot 2$. This gives $b=384\cdot 2$, exactly what we need. This yields $c=d=e=0$, $f=384\phi^{-1}$, $g=-384\phi$, $h=384\phi$, $i=384\cdot 1$, $j=-384\cdot 1$, $k=384\phi^{-1}$. Thus, if we start from a collection of the coefficients
\begin{eqnarray*}
(a,b,c,d,e,f,g,h,i,j,k)=(2,2,0,0,0,\phi^{-1},-\phi,\phi,1,-1,\phi^{-1}),
\end{eqnarray*} 
we arrive at the vector field $\frac{1}{384}(\varpi,\varrho,\sigma)$, where $(\varpi,\varrho,\sigma)$ is given by (\ref{pirmas}). Now, (\ref{one-lost}) gives\small
\begin{eqnarray*}
\mf{G}&=&2z\sin y+2y\sin z\\
&+&(-x+y-z)\sin\ell_{x,0}+(x-y-z)\sin\ell_{x,3}+(x+y+z)\sin\ell_{x,2}-(-x-y+z)\sin\ell_{x,1}\\
&+&\big{(}\phi x-y\big{)}
\sin\ell_{y,0}+\big{(}-\phi x+y\big{)}\sin\ell_{y,2}
+\big{(}-\phi x-y\big{)}\sin\ell_{y,1}-\big{(}\phi x+y\big{)}\sin\ell_{y,3}\\
&+&\big{(}-\phi^{-1}x+z\big{)}\sin\ell_{z,0}+\big{(}\phi^{-1}x+z\big{)}\sin\ell_{z,1}
+\big{(}\phi^{-1}x-z\big{)}\sin\ell_{z,3}-\big{(}-\phi^{-1}x-z\big{)}\sin\ell_{z,2}.
\end{eqnarray*}\normalsize
Finally, let us collect all functions in each row as factors of $x,y,z$, and use trigonometric addition formulas. For example,    
\begin{eqnarray*}
x\big{(}-\sin\ell_{x,0}+\sin\ell_{x,3}+\sin\ell_{x,2}+\sin\ell_{x,1}\big{)}&=&4x\sin\Big{(}\frac{x}{2}\Big{)}\sin\Big{(}\frac{\phi y}{2}\Big{)}\sin\Big{(}\frac{z}{2\phi}\Big{)}.
\end{eqnarray*}
Thus we get (in fact, MAPLE does these tedious computations for us) the first displayed formula in Theorem \ref{thm1}, where  $(\mf{V}_{x},\mf{V}_{y},\mf{V}_{z})=\frac{1}{2}(\mf{G}(x,y,z),\mf{G}(y,z,x),\mf{G}(z,x,y))$. For the first coordinate of its curl, we have $\mf{W}_{x}=\frac{\p\mf{V}_{z}}{\p y}-\frac{\p\mf{V}_{y}}{\p z}$, and this gives the second displayed formula in Theorem \ref{thm1}. To get immediately formulas in Theorem \ref{thm1}, in \cite{alkauskas-codes} we should set
\begin{eqnarray*}
(a,b,c,d,e,f,g,h,i,j,k)=\Big{(}1,1,0,0,0,\frac{\phi^{-1}}{2},-\frac{\phi}{2},\frac{\phi}{2},\frac{1}{2},-\frac{1}{2},\frac{\phi^{-1}}{2}\Big{)}.
\end{eqnarray*} 

  Some explanation is needed why the icosahedral symmetry holds for $\mf{W}$, too. This is clear from a coordinate-free definition of a curl; see, for example, (\cite{ficht3}, Chapter XVIII, \S 4). Indeed, let $F=(F_{x},F_{y},F_{z})$ be a smooth vector field, let $M$ be any point, and let $\m{n}$ be any direction from this point. In the plane perpendicular to it and passing through $M$, let us round the point $M$ with a region $\Sigma$, whose smooth boundary is $\lambda$ and an area is $|\Sigma|$. Then the Kelvin-Stokes theorem states that
\begin{eqnarray*}
(\nabla\times F)_{\m{n}}=\lim\limits_{\Sigma\rightarrow M}\frac{1}{|\Sigma|}\int\limits_{\lambda}F_{\lambda}\d\lambda,\quad F_{\lambda}\d\lambda=F_{x}\d x+F_{y}\d y+F_{z}\d z;
\end{eqnarray*}
here $(\nabla\times F)_{\m{n}}$ is a projection of a vector $\nabla\times F(M)$ onto a direction of $\m{n}$, and $\Sigma$ shrinks to the point $M$. Thus, we can define this projection coordinate-free, and since $\m{n}$ is any direction, this defines the curl. Now it is clear that if a vector field remains unchanged under a certain orthogonal transformation $\eta\in SO(3)$, so does its curl.\\

If we instead specialize
\begin{eqnarray} 
(a,b,c,d,e,f,g,h,i,j,k)=\Big{(}-2\phi^{-1},2\phi,0,0,0,-\phi^{-2},-\phi^{2},
-1,\phi,\phi^{-1},1\Big{)},
\label{Y}
\end{eqnarray}
we obtain Beltrami vector field $\mf{Y}$ with icosahedral symmetry whose Taylor series starts at degree $10$ rather than $6$, and $\mf{I}+a\mf{Y}$ is our $1$-parameter family.\\

Now we can start from the correct choice of parameters $a$ through $k$, and verify all the needed properties by MAPLE. The codes are provided by \cite{alkauskas-codes}. Changing parameters to the ones given by (\ref{Y}) verifies formulas and claims in Theorem 2.

\section{Final remarks}
\label{paskut}
Of course, such questions as
\begin{itemize}
\item[i)]which of the orbits under the flow with the vector field $\mf{I}$ are bounded; 
\item[ii)] are closed;
\end{itemize} 
 need to be investigated. Also,
\begin{itemize} 
\item[iii)]for which points $\m{x}\in\mathbb{R}^{3}$ there exists a bounded increasing sequence $\{t_{i}:i\in\mathbb{R}\}$, such that $F(\m{x},t_{i})$ tends to infinity;
\end{itemize}  
here $F(\m{x},t)$ is the flow with a vector field $\mf{I}$. All we can say that these classifications must have an icosahedral symmetry, too. \\

The existence of $C^{\infty}$ or even $C^{\omega}$ vector fields on a $3$-sphere $\m{S}^{3}$ with no circular orbits where demonstrated by Kuperberg \cite{kuperberg}, while for Beltrami $C^{\omega}$ vector fields (with non-vanishing curl) this was ruled out in \cite{etnyre}. It is interesting to know the answer for the flow $F(\m{x},t)$. At this moment, we can answer partially to these questions as follows.\\

\begin{figure}
\epsfig{file=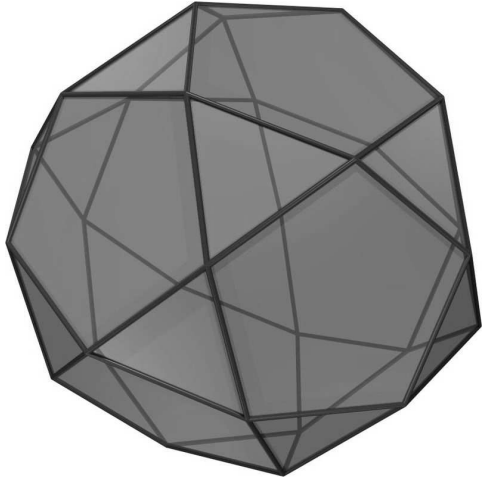,width=180pt,height=180pt,angle=0}
\caption{Icosidodecahedron}
\label{ico-dodeca}
\end{figure}

Consider the vector field $\mf{I}$. For the vector field $\mf{Y}$ the results are completely analogous, though the function $\Upsilon$ - see (\ref{ee}) - then is more complicated. Any point $\m{x}\in\mathbb{R}^{3}\setminus\{\m{0}\}$ has $60$ different points under the action of the group $\mathbb{I}$. However, there are three exceptional cases. Let $s\in\mathbb{R}\setminus\{0\}$.
\begin{itemize}
\item[\textbf{F}]- {\sc Faces}. The point $(\phi s,s,0)$ has 12 equivalent points under the action of $\mathbb{I}$;
\item[\textbf{V}]- {\sc Vertices}. The point $(s,s,s)$ has $20$ equivalent points;
\item[\textbf{E}]- {\sc Edges}. the point $(s,0,0)$ has $30$ equivalent points.
\end{itemize}
Geometrically, consider $6$ planes given by $(\phi^2x^2-y^2)(\phi^2 y^2-z^2)(\phi^2z^2-x^2)=0$. These $6$ planes split the unit sphere into $12$ pentagons having arcs of great circles as their sides, their centers being given by $\m{F}$ (more precisely, intersection of these lines with the unit sphere), $20$ triangles with centers being $\m{V}$, and $30$ intersection points given by $\m{E}$. Thus, it is more convenient to describe all these lines using not the geometry of a dodecahedron or an icosahedron, but rather an \emph{icosidodecahedron}; see Figure \ref{ico-dodeca}.\\

Each of $\m{F},\m{V},\m{E}$ (and its equivalent) defines a line passing through the origin, $62$ lines in total. On these $62$ lines the vector field $\m{M}$, as given by (\ref{pirmas}), vanishes. Each of these lines is divided into segments. Each segment is the complete orbit of the flow $F(\m{x},t)$.  Every orbit starts and finishes at two fixed points of the vector field $\mf{I}$, respectively. Since calculations are analogous in all three cases, we will show this in case \textbf{F}.\\

Indeed, then a point $(\phi s,s,0)$ has only $12$ equivalent points $(\pm\phi s,\pm s,0)$ (signs are independent), and all cyclic permutations. We take $(\phi s,s,0)$ as a representative. By a direct calculation, 
\begin{eqnarray*}
\mf{I}(\phi s,s,0)=\big{(}\phi\Upsilon(s),\Upsilon(s),0\big{)},
\end{eqnarray*}
where
\begin{eqnarray}
\Upsilon(s)=-s\sqrt{5}\Big{(}1-\phi\cos(s)+\phi^{-1}\cos(\phi s)\Big{)}. 
\label{ee}
\end{eqnarray}  
\begin{figure}
\centering
\begin{tabular}{c c}
\epsfig{file=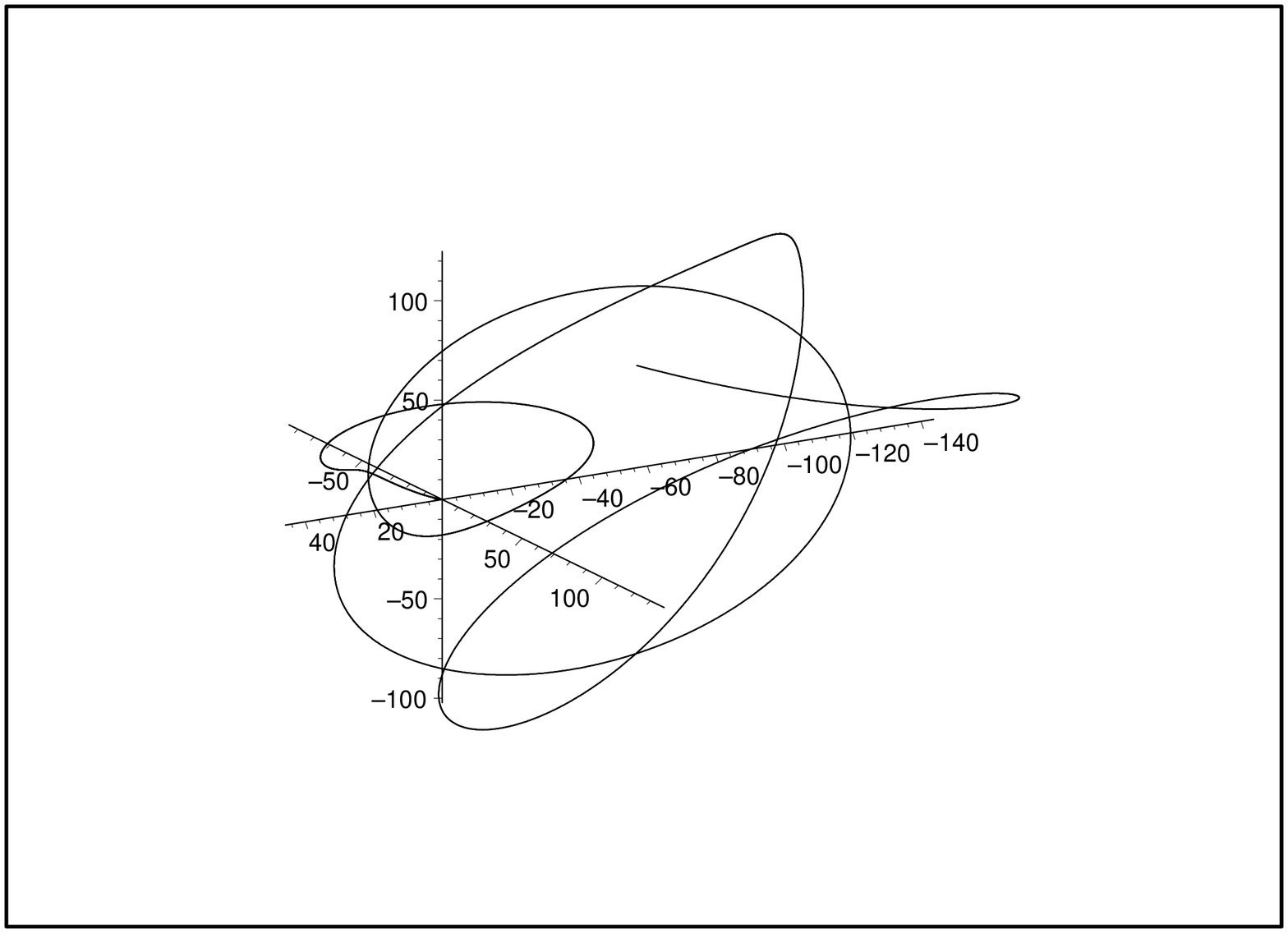,width=260pt,height=228pt,angle=-90}
&\epsfig{file=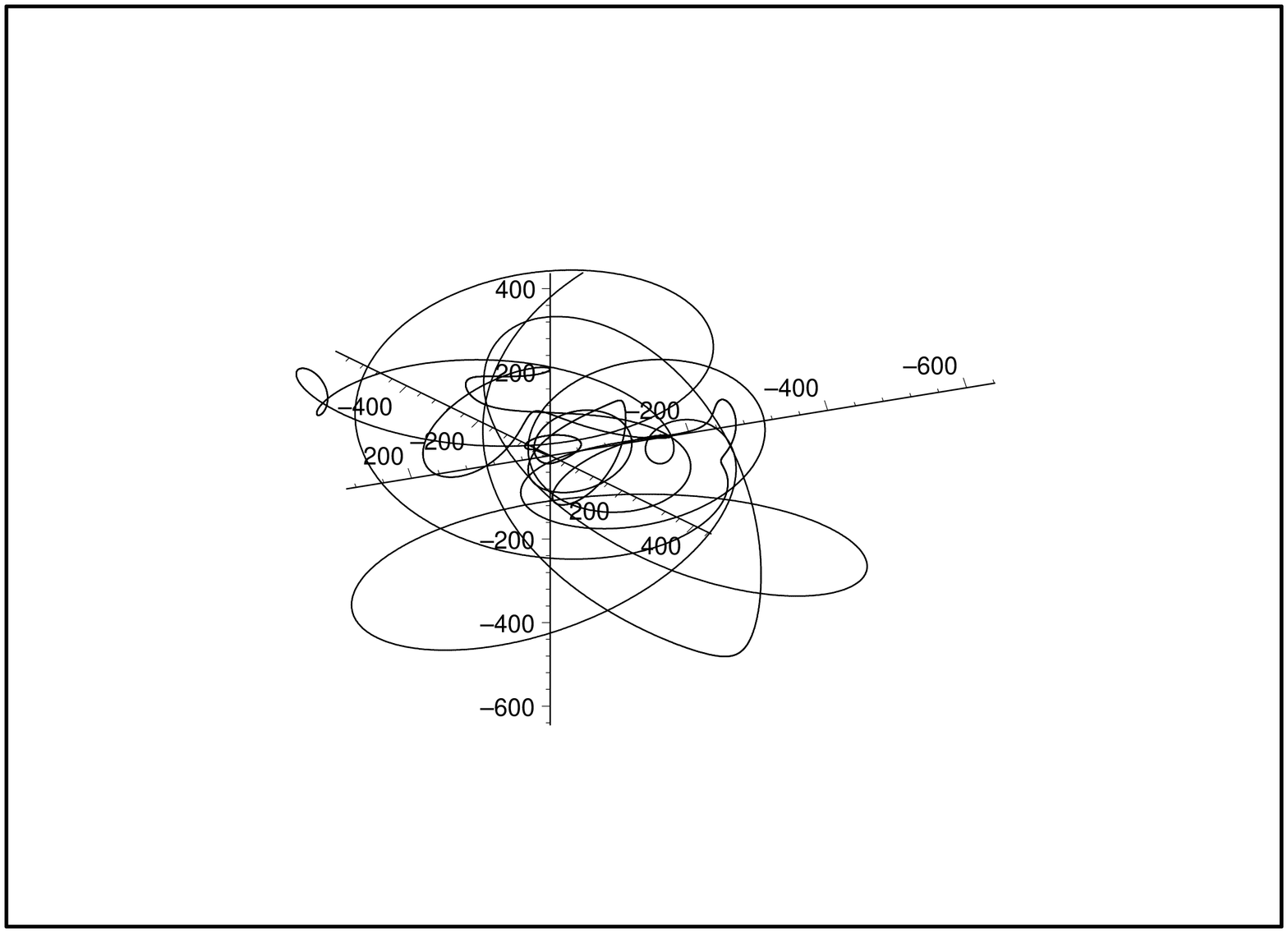,width=260pt,height=228pt,angle=-90} \\
The curve $\mf{I}(5s,6s,7s)$, $s\in[0,5]$.& The curve $\mf{I}(5s,6s,7s)$, $s\in[0,15]$.\\
\epsfig{file=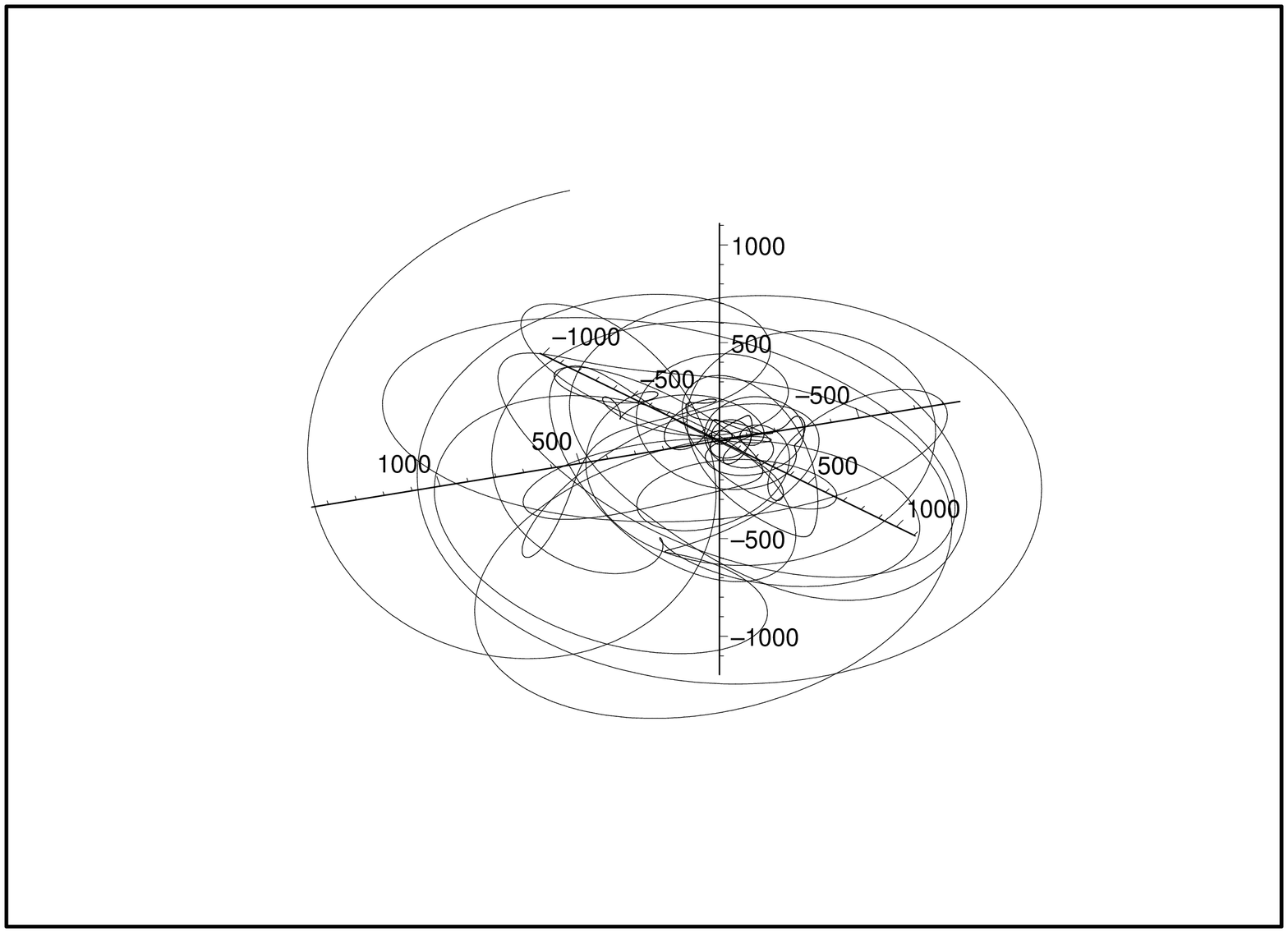,width=260pt,height=228pt,angle=-90} &
\epsfig{file=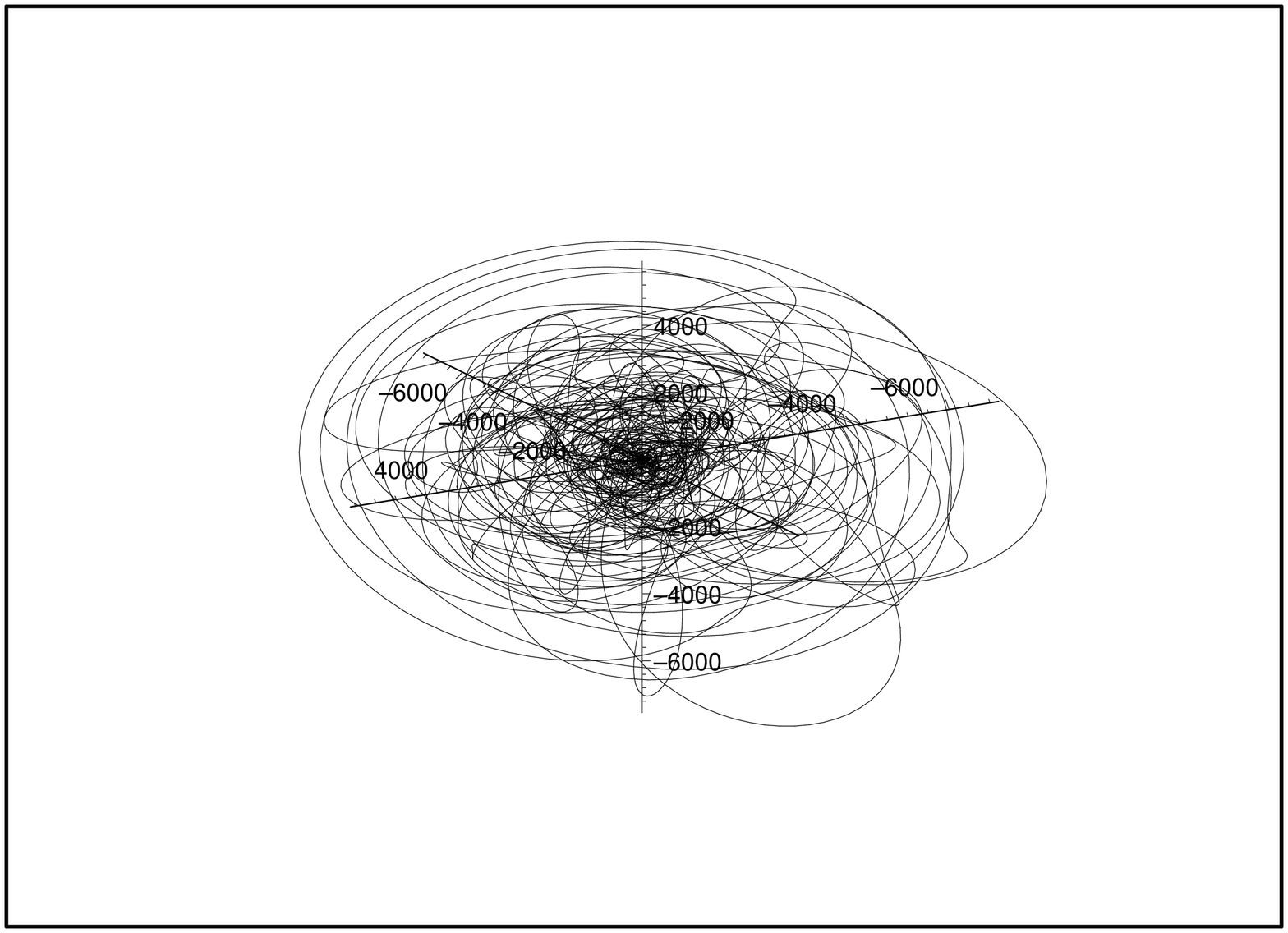,width=260pt,height=228pt,angle=-90}\\
The curve $\mf{I}(5s,6s,7s)$, $s\in[0,30]$. & The curve $\mf{I}(5s,6s,7s)$, $s\in[0,150]$.
\end{tabular}
\caption{Vector field $\mf{I}$ on a line $\mathbb{R}(5,6,7)$}.
\label{line}
\end{figure} 
This gives
\begin{eqnarray*}
\limsup\limits_{s\rightarrow\infty}\frac{\Upsilon(s)}{s}=2\sqrt{5}\phi^{-1},\quad\liminf\limits_{s\rightarrow\infty}\frac{\Upsilon(s)}{s}=-2\sqrt{5}\phi.
\end{eqnarray*}
Indeed, let $s=2\pi\ell$, $\ell\in\mathbb{N}$. Then $\phi s$ will hit arbitrarily close to $\pi+2\pi k$, $k\in\mathbb{N}$. For this we should have
\begin{eqnarray*}
2\pi\ell\phi=\pi+2\pi k+\epsilon(\ell)\Longrightarrow
\phi=\frac{1+2k}{2\ell}+\frac{\epsilon(\ell)}{2\pi\ell}.
\end{eqnarray*}
This will occur if $(2\ell,2k+1)$ is chosen to be a pair of two consecutive Fibonacci numbers $(F_{3n},F_{3n+1})$, $n\in\mathbb{N}$. From Diophantine properties of quadratic irrationals we know that then $\epsilon(\ell)\sim\frac{c}{\ell}$. This gives the $\limsup$ part, and analogously for $\liminf$.\\

 Thus we see that if $s_{1}$ and $s_{2}$ are two consecutive zeros of $\Upsilon(s)=0$, the segment joining $s_{1}\cdot(\phi,1,0)$ and $s_{2}\cdot(\phi,1,0)$ is the full orbit, and $\mf{I}(s_{1}(\phi,1,0))=\m{0}$, $F(s_{1}(\phi,1,0),t)=s_{1}(\phi,1,0)$. In this case the differential system (\ref{system}) turns out to be essentially
\begin{eqnarray*}
\dot{y}(t)=\Upsilon(y(t)).
\end{eqnarray*}
The behaviour of the vector field $\mf{I}$ on other lines passing through the origin is far more complicated. For example, Figure \ref{line} plots $\{\mf{I}(5s,6s,7s)$: $s\in[0,150]\}$; the choice of the line is motivated by a wish to further elucidate Figure \ref{fig-3}.\\

As a final remark of this paper, consider the Riemannian manifold $(\m{S}^{3},g)$ (a standard $3$-sphere), where $g$ is the usual induced Euclidean metric from $\mathbb{R}^{4}$, and curl operator was defined in the beginning of this paper. Now, consider a finite subgroup of $SO(4)$. For example, let $\mathbb{O}_{4}$ the the group of order $4!\cdot 2^{3}=192$, the so called \emph{orientation preserving hyperoctahedral group}, generated by matrices 
\begin{eqnarray*}
\alpha\mapsto\begin{pmatrix}
0 & 1 & 0 & 0\\
0 & 0 & 1 & 0\\
1 & 0 & 0 & 0\\
0 & 0 & 0 & 1\\
\end{pmatrix},
\beta\mapsto\begin{pmatrix}
1 & 0 & 0 & 0 \\
0 & -1 & 0 & 0 \\
0 & 0 & 0 & 1 \\
0 & 0 & 1 & 0 
\end{pmatrix},
\gamma\mapsto\begin{pmatrix}
0 & 1 & 0 & 0\\
1 & 0 & 0 & 0\\
0 & 0 & -1 & 0\\
0 & 0 & 0 & 1 \\
\end{pmatrix}.
\end{eqnarray*} 
We may ask for a similar question of constructing Beltrami vector field, equal to its own curl, with a $\mathbb{O}_{4}$-symmetry. This is our next task in a beautiful subject of constructing Beltrami vector fields with various symmetries.

\end{document}